\documentclass[12pt]{amsart}

\usepackage{amsfonts,amsmath,amssymb,amsthm}
\usepackage{latexsym}
\usepackage{mathrsfs}
\usepackage[hidelinks]{hyperref}
\usepackage{bookmark}

\newcommand{\gsl}{\mathfrak{sl}}

\newcommand{\osp}{\mathfrak{osp}}

\newcommand{\qi}{{q^{-1}}}
\newcommand{\Usl}{U_q(\gsl_2)}
\newcommand{\Utrp}{U_q(\gsl_2)^{\otimes 3}}

\newcommand{\cC}{\mathcal{C}}
\newcommand{\cR}{\mathcal{R}}
\newcommand{\cJ}{\mathcal{J}}
\newcommand{\bI}{\mathbb{I}}

\newcommand{\id}{\text{id}}
\newcommand{\End}{\text{End}}

\newcommand{\al}{\alpha}
\newcommand{\be}{\beta}
\newcommand{\ga}{\gamma}
\newcommand{\de}{\delta}

\newcommand{\C}{C}
\newcommand{\tX}{\widetilde{X}}
\newcommand{\tS}{\widetilde{S}}
\newtheorem{prop}{Proposition}[section]
\newtheorem{thm}{Theorem}[section]
\newtheorem{rem}{Remark}[section]
\newtheorem{corollary}{Corollary}[section]

\numberwithin{equation}{section}

\begin{document}

\title[]{Braid group and $q$-Racah polynomials}

\author[N.Cramp\'e]{Nicolas Cramp\'e$^{\dagger}$}
\address{$^\dagger$ Institut Denis-Poisson CNRS/UMR 7013 - Universit\'e de Tours - Universit\'e d'Orl\'eans, Parc de Grandmont, 37200 Tours, France.}
\email{crampe1977@gmail.com}

\author[L.Vinet]{Luc Vinet$^{*}$}
\author[M.Zaimi]{Meri Zaimi$^{*}$}
\address{$^*$ Centre de recherches math\'ematiques, Universit\'e de Montr\'eal, P.O. Box 6128, Centre-ville Station, Montr\'eal (Qu\'ebec), H3C 3J7, Canada.}
\email{vinet@crm.umontreal.ca, meri.zaimi@umontreal.ca}

\begin{abstract}
	The irreducible representations of two intermediate Casimir elements associated to the recoupling of three identical irreducible representations of $U_q(\mathfrak{sl}_2)$ are considered. It is shown that these intermediate Casimirs are related by a conjugation involving braid group representations. Consequently, the entries of the braid group matrices are explicitly given in terms of the $q$-Racah polynomials which appear as $6j$-symbols in the Racah problem for $U_q(\mathfrak{sl}_2)$. Formulas for these polynomials are derived from the algebraic relations satisfied by the braid group representations.      
\end{abstract}

\maketitle

\section{Introduction}\label{sec:intro}

This paper connects the representations of the braid group on three strands to the $q$-Racah polynomials by using the structure and the representation theory of the quantum algebra $\Usl$. The algebraic framework provided by the braid group allows to recover the orthogonality relation of these polynomials and to obtain in particular an addition formula.

The generic braid group is generated by $n-1$ invertible elements $s_i$ which commute for non-adjacent indices and which satisfy the braid relation $s_is_{i+1}s_i=s_{i+1}s_is_{i+1}$ otherwise \cite{A}. As its name suggests, this group has a diagrammatic representation in terms of the braiding of $n$ strands. Several algebras in mathematical physics can be seen as quotients of the braid group algebra; some examples are the symmetric group algebra, the Brauer, Hecke, Temperley--Lieb and Birman--Murakami--Wenzl algebras. The braid group representations play a role in the search for solutions to the Yang--Baxter equation and are related to the $R$-matrix \cite{J}.

The $q$-Racah polynomials form a finite family of basic hypergeometric orthogonal polynomials which sits at the top of the discrete $q$-Askey scheme \cite{Koek}. They can be obtained from a truncation of the Askey--Wilson polynomials. The bispectral properties of both families of polynomials are encoded in the Askey--Wilson algebra \cite{Zh}. An important feature of the $q$-Racah polynomials is that they appear as $6j$-symbols in the recoupling of three spin representations of the quantum algebra $\Usl$ \cite{GZ,KR}. These $6j$-symbols are the overlap coefficients between the bases which diagonalize the intermediate Casimirs associated to the coupling of the first two and last two factors of $\Utrp$. This is related to the realization of the Askey--Wilson algebra as the centralizer of the diagonal embedding of $\Usl$ in $\Utrp$ \cite{GZ,H}. 

The connection between the braid group and the $q$-Racah polynomials is established here by considering the recoupling of three identical spin representations of $\Usl$. On the one hand, the intermediate Casimirs are related by a transition matrix with entries given by $q$-Racah polynomials. On the other hand, the quasi-triangular Hopf algebra structure of $\Usl$ allows to define braided $R$-matrices which satisfy the braid relation on three strands. The crucial point of this paper is to observe that the intermediate Casimirs are also related by a conjugation involving these braided $R$-matrices. This observation is directly based on the results of a previous paper \cite{CGVZ} in which the formalism of the universal $R$-matrix is used to provide an intrinsic algebraic description of the intermediate Casimirs of $\Utrp$. The comparison between the transition matrix and the conjugation by braided $R$-matrices allows to connect the braid group to the $q$-Racah polynomials. This leads to a finite irreducible representation of the braid group which has appeared previously in the study and classification of spin Leonard pairs \cite{Curtin}.                 

The paper is structured as follows. Section \ref{sec:Uqsl2rep} presents some known facts on $\Usl$ and its representation theory. Section \ref{sec:braidqRacah} examines how the intermediate Casimirs of the tensor product of three identical irreducible $\Usl$-representations are conjugated in the two ways described above, and relates precisely the braid group representations to the $q$-Racah polynomials. Finally, Section \ref{sec:formulas} makes use of the algebraic interpretation of the $q$-Racah polynomials in terms of the braid group to obtain formulas involving these polynomials. The paper closes with concluding remarks and is complemented by some appendices.

\section{$U_q(\mathfrak{sl}_2)$ algebra and its representations}\label{sec:Uqsl2rep}

In this section, we briefly recall some fundamental aspects of the quantum algebra $\Usl$ and of its finite irreducible representations. Throughout this paper, we restrain to the case where $0<q<1$ to follow the conventions of \cite{Koek}.   

\subsection{$U_q(\mathfrak{sl}_2)$ algebra}\label{subsec:Uqsl2}

The associative algebra $\Usl$ is generated by $E$, $F$ and $q^{H}$ with the defining relations
\begin{equation}
	q^{H}E=q  Eq^H, \quad q^{H}F=q^{-1} Fq^H, \quad [E,F]=[2H]_q, \label{eq:Uqsl2}
\end{equation}
where $[X]_q=\frac{q^{X}-q^{-X}}{q-q^{-1}}$. The following Casimir element generates the center of $\Usl$
\begin{equation}
	\C=(q-\qi)^2 FE+qq^{2H}+\qi q^{-2H}. \label{eq:Casimir}
\end{equation}
There exists an algebra homomorphism $\Delta:\Usl \rightarrow \Usl \otimes \Usl $, called comultiplication, which is defined on the generators by 
\begin{equation}
	\Delta(E)=E \otimes q^{-H} +q^H \otimes E , \quad  \Delta(F)=F \otimes q^{-H} +q^{H} \otimes F, \quad \Delta(q^H)=q^H\otimes q^H, \label{eq:comult}
\end{equation}
and which is coassociative 
\begin{equation}
	(\Delta \otimes \id)\Delta= (\id \otimes \Delta)\Delta =: \Delta^{(2)}. \label{eq:coasso}
\end{equation}

\subsection{Finite irreducible representations of $U_q(\mathfrak{sl}_2)$ and tensor product decomposition}\label{subsec:irrepUqsl}

For each $j\in \{0,\frac{1}{2},1,\frac{3}{2},2,...\}$, the algebra $\Usl$ has an irreducible (spin-$j$) representation $M_j$ of finite dimension $2j+1$. We will denote the representation map by $\pi_{j}:\Usl \to \text{End}(M_j)$. The image by $\pi_{j}$ of the Casimir element $\C$ is proportional to the $(2j+1) \times (2j+1)$ identity matrix $\bI_{2j+1}$ : 
\begin{equation}
	\pi_{j}(\C)= \chi_j \bI_{2j+1}, \quad \text{where } \chi_j:= q^{2j+1}+q^{-2j-1}.
	\label{eq:casimirrep}
\end{equation}

Let us now fix an integer or half-integer $s$. In what follows, we will be interested in taking the tensor product of three copies of $M_s$. The tensor product of two irreducible spin-$s$ representations decomposes into the following direct sum of irreducible representations of $\Usl$
\begin{equation}
	M_s^{\otimes 2}=\bigoplus_{j=0}^{2s}M_j. \label{eq:decomp2}
\end{equation}
From this, one deduces the direct sum decomposition of the threefold tensor product
\begin{equation}
	M_s^{\otimes 3}=\bigoplus_{j\in \cJ_s}d_jM_j, \label{eq:decomp3}
\end{equation}
where $d_j\in\mathbb{Z}_{>0}$ is the degeneracy of $M_j$ and
\begin{equation}
	\cJ_s=\{j_{\min},j_{\min}+1,...,3s \}, \quad 
	j_{\min}=
	\begin{cases}
		0           & \text{if $s$ is integer,} \\
		\frac{1}{2} & \text{if $s$ is half-integer.} \\
	\end{cases} \label{eq:cJ}
\end{equation} 
The degeneracies $d_j$ can be computed explicitly:
\begin{equation}
	d_j=\min(2s,s+j) - |s-j| + 1 = 
	\begin{cases}
		2j+1  & \text{if } j_{\min} \leq j\leq s,\\
		3s-j+1 & \text{if } s<j\leq 3s.\\
	\end{cases} \label{eq:degen}
\end{equation}

\section{Braid group representations and $q$-Racah polynomials}\label{sec:braidqRacah}

In this section, we consider intermediate Casimirs of $\Usl$ in finite irreducible representations and show how they are related by braid group representations. This allows us to connect precisely those braid representations to the $q$-Racah polynomials.  

\subsection{Centralizer and intermediate Casimirs}\label{sec:Cent}

We define the centralizer $\cC_{s}$ of the image of the diagonal embedding of $\Usl$ in $\End(M_s^{\otimes 3})$
\begin{equation}
	\cC_{s}=\{m \in \End(M_s^{\otimes 3}) \  \big| \  [\pi_{s}^{\otimes 3} (\Delta^{(2)}(x)),m]=0,\forall x\in U_q(\gsl_2)\}. \label{eq:centr2}
\end{equation}
Using the properties of the comultiplication and the fact that $\C$ is central in $\Usl$, it is seen that the following elements belong to $\cC_s$
\begin{gather}
	C_1=\pi_{s}(\C)\otimes \bI_{(2s+1)^2}, \qquad C_2=\bI_{2s+1}\otimes \pi_{s}(\C) \otimes \bI_{2s+1}, \qquad C_3=\bI_{(2s+1)^2} \otimes \pi_{s}(\C), \label{eq:Ci} \\
	C_{12}=\pi_{s}^{\otimes 2}(\Delta(\C)) \otimes \bI_{2s+1}, \qquad C_{23}=\bI_{2s+1} \otimes \pi_{s}^{\otimes 2}(\Delta(\C)), \label{eq:Cij} \\
	C_{123}=\pi_{s}^{\otimes 3}(\Delta^{(2)}(\C)). \label{eq:C123}
\end{gather}
The elements $C_i$ and the total Casimir $C_{123}$ are central in $\cC_s$, while the intermediate Casimirs $C_{12}$ and $C_{23}$ satisfy the relations (of a version) of the Askey--Wilson algebra \cite{GZ}
\begin{align}
	\frac{[C_{23},[C_{12},C_{23}]_q]_q}{(q-\qi)^2}=&(q+\qi)^2C_{12}+(C_1 C_3 +C_2 C_{123})C_{23}-(q+\qi)(C_1C_2+C_3 C_{123}), \label{eq:AW1}\\
	\frac{[[C_{12},C_{23}]_q,C_{12}]_q}{(q-\qi)^2}=&(q+\qi)^2C_{23}+(C_1 C_3 +C_2 C_{123})C_{12}-(q+\qi)(C_2C_3+C_1 C_{123}), \label{eq:AW2}
\end{align}
where $[X,Y]_q=qXY-\qi YX$ is the $q$-commutator. Note that equations \eqref{eq:casimirrep} and \eqref{eq:Ci} imply that $C_1=C_2=C_3=\chi_s\mathbb{I}_{(2s+1)^3}$.

\subsection{Finite irreducible representations of the centralizer and $q$-Racah polynomials}\label{subsec:irrepCent}

The threefold tensor product of the space $M_s$ discussed in Subsection \ref{subsec:irrepUqsl} can also be decomposed as
\begin{equation}
	M_s^{\otimes 3}=\bigoplus_{j\in \cJ_s} M_j \otimes V_j, \label{eq:decompV}
\end{equation}
where the spaces $V_j$ for $j\in \cJ_s$ are all the non-equivalent finite irreducible representations of the centralizer $\cC_{s}$. Let us now fix some $\ell\in\cJ_s$. The space $V_\ell$ is of dimension $d_\ell$, which is given explicitly in \eqref{eq:degen}, and it can be constructed from the highest weight vectors of the multiple copies of the representation $M_\ell$ in the decomposition \eqref{eq:decomp3}. The representations of the elements $C_{12}$ and $C_{23}$ in $\End(V_\ell)$ can be found with the help of the Askey--Wilson algebra, as it was first shown in \cite{GZ}. We will only state here the results that are relevant for our purposes\footnote{For more details on how the computations can be reproduced for the specific case that we are considering here, see the proof of Proposition 4.2 in \cite{CVZ2}; the results of the present paper are obtained directly by taking $j_1=j_2=j_3=s$.}.

Let us denote for convenience
\begin{equation}
	a = |s-\ell|, \quad  N=d_\ell-1=\text{min}( 2s , s+\ell ) - |s-\ell|. \label{eq:aN}
\end{equation}
In the irreducible representation $V_\ell$, the element $C_{i}$ for $i\in\{1,2,3\}$ (resp. $C_{123}$) acts as the constant matrix of value $\chi_s$ (resp. $\chi_\ell$). We can choose a basis $\{v_j\}_{j=0}^{N}$ of $V_\ell$ such that the intermediate Casimirs $C_{12}$ and $C_{23}$ are represented respectively by the following matrices $X_1$ and $X_2$: 
\begin{equation}
	X_1=
	\begin{pmatrix}
		\chi_{a} &       &       &        &         \\
		& \chi_{a+1} &       &        &         \\
		&       & \chi_{a+2} &        &         \\
		&       &       & \ddots &         \\
		&       &       &        & \chi_{a+N}
	\end{pmatrix}, \quad 
	X_2=
	\begin{pmatrix}
		\al_{0,0} & \al_{0,1} &            &             &             \\
		\al_{1,0} & \al_{1,1} & \al_{1,2}  &             &             \\
		& \al_{2,1} & \al_{2,2}  & \ddots      &             \\
		&           & \ddots     & \ddots      & \al_{N-1,N} \\
		&           &            & \al_{N,N-1} & \al_{N,N}
	\end{pmatrix}, \label{eq:Xi}
\end{equation}
where
\begin{align}
	&\al_{j-1,j}=q^{-2a-1}\frac{(1-q^{2(j+a)})(1-q^{2(j+2a)})(1-q^{2(j-N-1)})(1- q^{2(j+3a+N+1)})}{(1-q^{2(2j+2a-1)})(1-q^{2(2j+2a)})}, \label{eq:aljmj}\\
	&\al_{j,j-1}=q^{2a+1}\frac{(1-q^{2(j+a)})(1-q^{2j})(1-q^{2(j+2a+N+1)})(1-q^{2(j-a-N-1)})}{(1-q^{2(2j+2a)})(1-q^{2(2j+2a+1)})}, \label{eq:aljjm} \\
	&\al_{j,j}=\chi_a-(\al_{j,j+1}+\al_{j,j-1}). \label{eq:aljj}
\end{align}
 
It will be useful to define the diagonal matrix
\begin{equation}
	D:=\text{diag}(\be_0,\be_1,...,\be_N) \label{eq:D}
\end{equation}
with
\begin{equation}
	\be_{i}:=\prod_{k=1}^{i}\sqrt{\frac{\al_{k-1,k}}{\al_{k,k-1}}}, \quad i=0,1,...,N. \label{eq:betai}
\end{equation}
One then verifies that the following matrix is tridiagonal and symmetric:
\begin{equation}
	\tX_2:=DX_2D^{-1}. \label{eq:X2t}
\end{equation}

Both $X_1$ and $X_2$ have the same eigenvalues $\chi_{a+j}$ for $j=0,1,...,N$. As it is known, these finite irreducible representations of the intermediate Casimirs $C_{12}$ and $C_{23}$ are related by a change of basis which involves the $6j$-symbols for $\Usl$. The next proposition gives precisely this change of basis for the specific case that concerns us. Note that we will use the standard notations for the $q$-shifted factorials, the $q$-hypergeometric functions ${}_r\phi_s$ and the $q$-Racah polynomials \cite{GR,Koek}. We have included the definitions in Appendix \ref{sec:appendix1}, and we have also listed there some useful identities.  

\begin{prop}\label{prop:diagX2} \cite{KR}
	In the irreducible representation $V_\ell$ of the centralizer $\cC_s$, the similarity transformation which relates the tridiagonal matrix $X_2$ to the diagonal matrix $X_1$ is
	\begin{equation}
		X_2=PX_1P^{-1}, \label{eq:diagX2}
	\end{equation}
	where $P$ is the invertible matrix whose entries are given by the following renormalized $q$-Racah polynomials
	\begin{align}
		P_{ij}&=\be_j^2 \, R_i(\mu(j);q^{2a},q^{2a},q^{2(3a+N+1)},q^{-2(a+N+1)}|q^2) \quad \text{for } i,j=0,1,...,N, \label{eq:qRacparam} \\
		&=\be_j^2 \, {}_{4}\phi_{3}\biggl(\genfrac..{0pt}{}{q^{-2i},q^{2(2a+1+i)},q^{-2j},q^{2(2a+1+j)}}{q^{2(a+1)},q^{-2N},q^{2(3a+N+2)}};q^2,q^2\biggr),  \label{eq:compoP}
	\end{align}
	and where the parameters $a$ and $N$ are given in terms of the fixed spins $s$ and $\ell$ in \eqref{eq:aN}. The normalization is explicitly given by
	\begin{equation}
		\beta_j^2=q^{-2j(2a+1)}\frac{(q^{2a+3},-q^{2a+3},q^{2(2a+1)},q^{2(3a+N+2)},q^{-2N};q^2)_j}{(q^2,q^{2a+1},-q^{2a+1},q^{-2(a+N)},q^{2(2a+N+2)};q^2)_j}. \label{eq:betar2}
	\end{equation}
\end{prop}
\proof 
For some fixed $j\in \{0,1,...,N\}$, let $u(j)=\sum_{i=0}^{N} u_i(j)v_i$ be the eigenvector of $X_2$ with eigenvalue $\chi_{a+j}$, where $v_i$ are the basis vectors of $V_\ell$ defined in Subsection \ref{subsec:irrepCent}. By construction, the matrix $P$ with entries $P_{ij}=u_i(j)$ is such that relation \eqref{eq:diagX2} holds. The eigenvalue equation $X_2 u(j)=\chi_{a+j} u(j)$ leads to the following three-term recurrence relation
\begin{equation}
	\al_{i,i-1} u_{i-1}(j) +\al_{i,i} u_i(j)+\al_{i,i+1}u_{i+1}(j)=\chi_{a+j} u_i(j), \quad i=0,1,...,N, \label{eq:EV}
\end{equation}
where we have used \eqref{eq:Xi}, with the coefficients $\al_{i,j}$ defined in \eqref{eq:aljmj}--\eqref{eq:aljjm}. Note that these definitions imply the boundary conditions $\al_{0,-1}=\al_{N,N+1}=0$.
It is straightforward to verify that \eqref{eq:EV} is the recurrence relation for the $q$-Racah polynomials given in \eqref{eq:qRacparam} (see \cite{Koek}). The renormalization by $\be_j^2$ will be convenient. The explicit expression \eqref{eq:betar2} is obtained by using the definitions \eqref{eq:betai} and \eqref{eq:aljmj}--\eqref{eq:aljjm}, and  the identity \eqref{eq:qfid4}.
\endproof

\subsection{Braid group representations} 
The matrices $X_i$ satisfy the following minimal polynomial equation
\begin{equation}
	\prod_{j=0}^{N}(X_i-\chi_{a+j}\bI_{N+1})=0 \quad \text{for } i=1,2. \label{eq:minpolX}
\end{equation}
This relation allows to define, for $i=1,2$ and $r= 0,1,...,N$, the orthogonal idempotents
\begin{equation}
	E_i^{(r)} := \prod_{\genfrac{}{}{0pt}{}{k=0}{k\neq r}}^{N} \frac{X_i-\chi_{a+k}\bI_{N+1}}{\chi_{a+r}-\chi_{a+k}} \label{eq:idemp}
\end{equation}
which satisfy 
\begin{equation}
	E_i^{(r)}E_i^{(p)}=\delta_{rp}E_i^{(r)}, \quad \sum_{r=0}^{N} E_i^{(r)}=\bI_{N+1},  \quad E_i^{(r)}X_i=X_iE_i^{(r)}=\chi_{a+r} E_i^{(r)}. \label{eq:orthcompl}
\end{equation}
We also define the matrices $S_i\in \End(V_\ell)$, for $i=1,2$, by using these idempotents
\begin{equation}
	S_i:=\sum_{r=0}^{N}  \ga_r  E_i^{(r)}, \quad \text{where }\ga_r:=(-1)^{r} q^{r(r+2a+1)}. \label{eq:Si}
\end{equation}
The matrices $S_i$ are invertible and their inverses are given by $S_i^{-1}=\sum_{r=0}^{N}  \ga_r^{-1}  E_i^{(r)}$. 
From the definitions \eqref{eq:idemp} and \eqref{eq:Si}, it is seen that $S_i$ is a polynomial in $X_i$
\begin{equation}
	S_i=\sum_{r=0}^{N}  \ga_r  \prod_{\genfrac{}{}{0pt}{}{k=0}{k\neq r}}^{N} \frac{X_i-\chi_{a+k}\bI_{N+1}}{\chi_{a+r}-\chi_{a+k}}. \label{eq:SipolXi}
\end{equation}
As a consequence, $S_1$ is a diagonal matrix with diagonal entries $\ga_r$ for $r=0,1,...,N$. 

The algebra $\Usl$ is a quasi-triangular Hopf algebra and hence possesses a universal $R$-matrix which satisfies in particular the Yang--Baxter equation (see Appendix \ref{ssec:Runi}). In \cite{CGVZ}, it is shown that the intermediate Casimir elements of $\Utrp$ can be described with the help of the universal $R$-matrix of $\Usl$. In the tensor product of three identical irreducible $\Usl$-representations, this universal $R$-matrix allows to define two braided $R$-matrices $\check{R}_i$ which are closely related to the matrices $S_i$ defined by \eqref{eq:Si} and which satisfy the braid relation in the form of the braided Yang--Baxter equation (see Appendix \ref{ssec:braidR}). On the basis of these remarks, the following proposition gives additional algebraic relations satisfied by the matrices $S_i$ and the representations $X_i$ of the intermediate Casimirs $C_{12}$ and $C_{23}$ that will play a key role in the remaining of this paper. 
\begin{prop}\label{prop:relSiXi}
	The matrices $S_i$ satisfy the braid relation
	\begin{equation}
		S_1S_2S_1=S_2S_1S_2. \label{eq:Sibraid}
	\end{equation}
	Moreover, the tridiagonal matrix $X_2$ is related to the diagonal matrix $X_1$ as follows 
	\begin{equation}
		X_2= S_{1}S_{2}  X_1 S_{2}^{-1}S_{1}^{-1}=S_1^{-1}S_{2}^{-1}  X_1 S_{2}S_{1}. \label{eq:X12Si}
	\end{equation}
\end{prop}
\proof
The idempotent $E_1^{(r)}$ (resp. $E_2^{(r)}$) is the representation in $\End(V_\ell)$ of the element in $\cC_s$ which acts trivially on the third (resp. first) factor of $M_s^{\otimes 3}$, and which acts on the two other factors as the projector on the representation $M_{a+r}$ of the decomposition \eqref{eq:decomp2}. Hence, using the expression of the eigenvalue $\ga_r$ given in \eqref{eq:Si}, it can be shown that $S_i$ is proportional to a representation of the braided $R$-matrix $\check{R}_i$ of $\Usl$ \cite{LZ}. Relation \eqref{eq:Sibraid} then follows from the fact that the braided $R$-matrices satisfy the braided Yang--Baxter equation. For more details, see Appendix \ref{sec:appendix2}. The two equalities in \eqref{eq:X12Si} are a consequence of Theorem 3.1 in \cite{CGVZ}, which describes two intermediate Casimirs of $\Usl$ associated to the coupling of the first and third factors of $\Utrp$ as conjugations by a universal $R$-matrix; writing these conjugations in terms of braided $R$-matrices leads to relations between the intermediate Casimirs $C_{12}$ and $C_{23}$ which are represented by \eqref{eq:X12Si} in the irreducible representation $V_\ell$ of the centralizer $\cC_s$. 
\endproof

\subsection{Relating the braid group representations to the $q$-Racah polynomials}
We now combine the results of Propositions \ref{prop:diagX2} and \ref{prop:relSiXi} to get the main result of this section. 
\begin{thm}\label{thm:S2P}
	The braid matrices $S_i$ are related to the transition matrix $P$ as follows
	\begin{align}
		S_1S_2S_1&=m_qP, \label{eq:S2P} \\
		(S_1S_2S_1)^{-1}&=m_{\qi}P,     \label{eq:S2iP}
	\end{align}
	where
	\begin{align}
		m_q&=q^{2N(3a+N+2)}\frac{(q^{-2(a+N)};q^2)_N}{(q^{4(a+1)};q^2)_N}. \label{eq:m}
	\end{align}
\end{thm}
\proof
Relation \eqref{eq:diagX2} can be written in terms of the symmetric matrix $\tX_2$ \eqref{eq:X2t} and the diagonal matrix $D$ \eqref{eq:D} as follows
\begin{equation}
	\tX_2=(DPD^{-1})X_1(DPD^{-1})^{-1}. \label{eq:diagX2t}
\end{equation}
We also define $\tS_2:=DS_2D^{-1}$, which is symmetric by \eqref{eq:SipolXi}, and we write the first equality of \eqref{eq:X12Si} as
\begin{equation}
	\tX_2= (S_{1}\tS_2S_1)  X_1 ( S_{1}\tS_2S_1)^{-1}, \label{eq:X12Sit}
\end{equation}
where we have used the fact that $S_1X_1=X_1S_1$. Comparing \eqref{eq:diagX2t} with \eqref{eq:X12Sit}, one deduces that 
\begin{equation}
	S_{1}\tS_2S_1=DPD^{-1}M, \label{eq:SiDPM}
\end{equation}
where $M$ is a matrix which commutes with $X_1$. Since $X_1$ is diagonal with distinct eigenvalues, $M$ must be diagonal. The L.H.S. of equation \eqref{eq:SiDPM} is symmetric and therefore we should have $(DPD^{-1}M)^T=DPD^{-1}M$. The matrix $DPD^{-1}$ has entries $\be_iP_{ij}\be_j^{-1}$ and is hence symmetric, as it is seen from \eqref{eq:compoP}. Since $D$ and $M$ are diagonal, one thus finds that $MP=PM$. This equation is only possible if $M$ is proportional to the identity by some constant $m_q$, as it can be shown by writing the equation in matrix components and using the fact that the components $P_{ij}$ are never zero. One then recovers \eqref{eq:S2P} from \eqref{eq:SiDPM}.

To determine the constant $m_q$, one can use \eqref{eq:S2P} in the braid relation \eqref{eq:Sibraid} and examine the $(0,0)$ matrix component of the resulting equation. One finds 
\begin{equation}
	m_q=\left(\sum_{r=0}^{N}\ga_r^{-1}\be_r^2 \right)^{-1}. \label{eq:constants}
\end{equation}
It is straightforward to show by using the expression \eqref{eq:betar2} for $\beta_r^2$, the definition of $\ga_r$ given in \eqref{eq:Si} and the definition \eqref{eq:rphis} of the $q$-hypergeometric series that
\begin{align}
	\sum_{r=0}^{N}\ga_r^{-1}\beta_r^2&= {}_{6}\phi_{4}\biggl(\genfrac..{0pt}{}{q^{2a+3},-q^{2a+3},q^{2(2a+1)},q^{2(3a+N+2)},0,q^{-2N}}{q^{2a+1},-q^{2a+1},q^{-2(a+N)},q^{2(2a+N+2)}};q^2,q^{-2(3a+2)}\biggr). \label{eq:dr64}
\end{align}
The following formula, given in \cite{Koek}, is obtained by taking a limit of the Jackson summation formula for a terminating very-well-poised ${}_{6}\phi_{5}$ series 
\begin{equation}
	{}_{6}\phi_{4}\biggl(\genfrac..{0pt}{}{qa_1^{\frac{1}{2}},-qa_1^{\frac{1}{2}},a_1,a_2,0,q^{-n}}{a_1^{\frac{1}{2}},-a_1^{\frac{1}{2}},qa_1a_2^{-1},a_1q^{n+1}};q,\frac{q^{n}}{a_2}\biggr)=\frac{(qa_1;q)_n}{(qa_1a_2^{-1};q)_n}a_2^{-n}, \quad n=0,1,2,... \label{eq:6phi4}
\end{equation}   
Equation \eqref{eq:6phi4} with the substitutions $q\to q^2$, $a_1=q^{2(2a+1)}$ and $a_2=q^{2(3a+N+2)}$ implies
\begin{equation}
	\sum_{r=0}^{N}\ga_r^{-1}\be_r^2
	=\frac{(q^{4(a+1)};q^2)_N}{(q^{-2(a+N)};q^2)_N}q^{-2N(3a+N+2)}. \label{eq:sumdr2}
\end{equation}
Result \eqref{eq:m} then directly follows from \eqref{eq:constants} and \eqref{eq:sumdr2}.

Finally, to prove \eqref{eq:S2iP}, one can use the fact that the matrix $P$ is invariant under $q\to \qi$. This can be shown by using expressions \eqref{eq:compoP} and \eqref{eq:betar2}, the definition of $q$-hypergeometric series given in \eqref{eq:rphis}, and the identity \eqref{eq:qfid2}. One can also show that $S_i \to S_i^{-1}$ when $q\to \qi$, for $i=1,2$. This is obvious for $S_1$, which is diagonal with eigenvalues $\ga_r$ given in \eqref{eq:Si}, while for $S_2$, this follows from the fact that $S_2=PS_1P^{-1}$, as seen from \eqref{eq:diagX2} and \eqref{eq:SipolXi}. Equation \eqref{eq:S2iP} is then obtained by taking $q\to \qi$ in \eqref{eq:S2P}. Alternatively, one could proceed in a similar manner as for the proof of \eqref{eq:S2P}, by using the second equality of \eqref{eq:X12Si} and the inverse of the braid relation \eqref{eq:Sibraid} instead. 
\endproof

\begin{corollary}\label{cor:repbraid}
	The diagonal matrix $S_1=\text{diag}(\ga_0,\ga_1,...,\ga_N)$ and the matrix $S_2$ with entries
	\begin{equation}
		(S_2)_{ij}=m_q\ga_i^{-1}\ga_j^{-1}\be_j^2 \, R_i(\mu(j);q^{2a},q^{2a},q^{2(3a+N+1)},q^{-2(a+N+1)}|q^2) \quad \text{for } i,j=0,1,...,N \label{eq:Sij}
	\end{equation}
	provide a representation of the braid group on three strands.
\end{corollary}
\proof
This is a direct consequence of the braid relation \eqref{eq:Sibraid} in Proposition \ref{prop:relSiXi}. The matrix entries given in \eqref{eq:Sij} follow from the equation \eqref{eq:S2P} in Theorem \ref{thm:S2P} and from the expression \eqref{eq:qRacparam} in Proposition \ref{prop:diagX2}. 
\endproof

\begin{rem} 
	With reference to the previous works in \cite{Curtin,NT}, the matrices $X_i$ form a Leonard pair, as it is seen from their definition \eqref{eq:Xi} and equation \eqref{eq:diagX2}. Moreover, the pair of invertible matrices $S_1,S_2$ is a Boltzmann pair for $X_1,X_2$; this is a consequence of definition \eqref{eq:Si} and the relation \eqref{eq:X12Si}. Therefore, $X_1,X_2$ is an example of a spin Leonard pair. More precisely, it is a spin Leonard pair of type I with the parameters of \cite{Curtin} changed as follows : $q \to q^2$, $d\to N$, $\theta_0\to\chi_a$, $h\to q^{-2a-1}$ and $\nu\to -q^{2a+2}$. The orthogonal polynomials appearing in \eqref{eq:Sij} correspond to those appearing in the description of Boltzmann pairs for spin Leonard pairs of type I in \cite{Curtin}. The braid relation \eqref{eq:Sibraid} also arises in the context of spin Leonard pairs and implies that $S_1,S_2$ is a balanced Boltzmann pair. The fact that spin Leonard pairs are self-dual is reflected here by the fact that the matrix $P^2$ is proportional to the identity, as will be discussed in the next section.      	
\end{rem}

\section{Formulas for $q$-Racah polynomials}\label{sec:formulas}

In this section, we use the results of Theorem \ref{thm:S2P} and the algebraic relations satisfied by the matrices $S_i$ and $X_i$ in order to obtain formulas involving $q$-Racah polynomials.
\subsection{Orthogonality relation} We can recover the orthogonality relation of the $q$-Racah polynomials from the matrix equation $S_2S_2^{-1}=\bI_{N+1}$. Using results \eqref{eq:S2P} and \eqref{eq:S2iP} to substitute $S_2$ and $S_2^{-1}$ in the previous equation, we obtain
\begin{equation}
	m_qm_{q^{-1}}P^2=\bI_{N+1}.
\end{equation}
Writing this equation in matrix components gives 
\begin{align}
	\sum_{k=0}^{N}P_{ik}P_{kj} & =  (m_qm_{q^{-1}})^{-1}\delta_{ij}. \label{eq:orthre1}
\end{align}
We now substitute expressions \eqref{eq:qRacparam}, \eqref{eq:betar2} and \eqref{eq:m} in \eqref{eq:orthre1}, and we take $q\to q^{1/2}$ to get
\begin{align}
	&\sum_{k=0}^{N}\frac{(1-q^{2a+1+2k})}{q^{k(2a+1)}(1-q^{2a+1})}\frac{(q^{2a+1},q^{3a+N+2},q^{-N};q)_k}{(q,q^{-(a+N)},q^{2a+N+2};q)_k}R_i(\mu(k))R_j(\mu(k)) \label{eq:orthre2} \\
	& =  \frac{(q^{2a+2},q^{-2a-N-1};q)_N}{(q^{-a-N},q^{a+1};q)_N}\frac{q^{i(2a+1)}(1-q^{2a+1})}{(1-q^{2a+1+2i})}\frac{(q,q^{-a-N},q^{2a+N+2};q)_i}{(q^{2a+1},q^{3a+N+2},q^{-N};q)_i}\delta_{ij}, \nonumber
\end{align}
where we have used the notation
\begin{equation}
	R_i(\mu(j))=R_i(\mu(j);q^{a},q^{a},q^{3a+N+1},q^{-a-N-1}|q). \label{eq:Rij}
\end{equation}
We have also used identities \eqref{eq:qfid2}, \eqref{eq:qfid3} and \eqref{eq:qfid5}, and the fact that $R_i(\mu(j))=R_j(\mu(i))$ to arrive at \eqref{eq:orthre2}. 
It is straightforward to verify that \eqref{eq:orthre2} is precisely the orthogonality relation of the $q$-Racah polymomials given in \eqref{eq:Rij} (see \cite{Koek}).

\subsection{Addition formula} We now use the braid relation satisfied by the matrices $S_i$ in order to obtain an addition formula for ${}_{4}\phi_{3}$ functions.
\begin{prop}\label{prop:addform}
	For $a,N=0,1,2,...$ and $i,j=0,1,...,N$, the following formula holds
	\begin{align}
		&\sum_{k=0}^N c_k\, {}_{4}\phi_{3}\biggl(\genfrac..{0pt}{}{q^{-i},q^{2a+i+1},q^{-k},q^{2a+k+1}}{q^{a+1},q^{-N},q^{3a+N+2}};q,q\biggr){}_{4}\phi_{3}\biggl(\genfrac..{0pt}{}{q^{-k},q^{2a+k+1},q^{-j},q^{2a+j+1}}{q^{a+1},q^{-N},q^{3a+N+2}};q,q\biggr) \label{eq:addform}\\
		=&(-1)^{i+j-N}q^{{i\choose 2}+{j\choose 2}-{N\choose 2}}q^{(i+j-2N)(a+1)}\frac{(q^{2a+2};q)_N}{(q^{a+1};q)_N} {}_{4}\phi_{3}\biggl(\genfrac..{0pt}{}{q^{-i},q^{2a+i+1},q^{-j},q^{2a+j+1}}{q^{a+1},q^{-N},q^{3a+N+2}};q,q\biggr), \nonumber
	\end{align}
	where
	\begin{equation}
		c_k=(-1)^{k}q^{-{k\choose 2}}q^{-k(3a+2)}\frac{(1-q^{2a+1+2k})}{(1-q^{2a+1})}\frac{(q^{2a+1},q^{3a+N+2},q^{-N};q)_k}{(q,q^{-a-N},q^{2a+N+2};q)_k}.
	\end{equation}
	In terms of $q$-Racah polynomials, this can be written as
	\begin{equation}
		\sum_{k=0}^N c_k\, R_i(\mu(k))R_j(\mu(k))=(-1)^{i+j-N}q^{{i\choose 2}+{j\choose 2}-{N\choose 2}}q^{(i+j-2N)(a+1)}\frac{(q^{2a+2};q)_N}{(q^{a+1};q)_N} R_i(\mu(j)),
	\end{equation}
	where we have used again the notation given in \eqref{eq:Rij}.
\end{prop}
\proof
One can use result \eqref{eq:S2P} to write the braid relation \eqref{eq:Sibraid} in terms of the matrix $P$ and get
\begin{equation}
	S_1PS_1=m_qPS_1^{-1}P.
\end{equation}
Writing the previous equation in matrix components, one finds
\begin{equation}
	\ga_i\ga_jP_{ij}=m_q\sum_{k=0}^N \ga_k^{-1}P_{ik}P_{kj}.
\end{equation}
It is then straightforward to obtain \eqref{eq:addform} by using the explicit expressions of $P_{ij}$, $\be_j^2$, $\ga_k$ and $m_q$, given respectively in \eqref{eq:compoP}, \eqref{eq:betar2}, \eqref{eq:Si} and \eqref{eq:m}. Note that we have used identity \eqref{eq:qfid3}, and we have taken $q\to q^{1/2}$ at the end of the computation to write the identity in a nicer form.
\endproof

\begin{rem}
	Formula \eqref{eq:addform} in Proposition \ref{prop:addform} has a structure similar to the q-analog of the Racah identity \cite{BL,Rac} obtained in \cite{KR}.
	%\footnote{Note that there is a typo in the formula (6.17) in \cite{KR} : the spins $j_2$ and $j_3$ should be exchanged in the second $6j$-symbol inside the summation.}.
	It did not prove possible to cast relation \eqref{eq:addform} as a special case of the latter nor of the addition formula for the $q$-Racah polynomials that results from the approach presented in \cite{Rosen}.
\end{rem}

\subsection{$q$-Racah polynomials in terms of recurrence coefficients} 
Let us finally remark that by using result \eqref{eq:S2P} to write $S_2$ in terms of $P$ in the relation \eqref{eq:SipolXi}, we find
\begin{equation}
	P=m_q^{-1}S_1\left( \sum_{r=0}^{N}  \ga_r  \prod_{\genfrac{}{}{0pt}{}{k=0}{k\neq r}}^{N} \frac{X_2-\chi_{a+k}\bI_{N+1}}{\chi_{a+r}-\chi_{a+k}}\right)S_1. \label{eq:PX2}
\end{equation}
Written in matrix components, this equation provides an expression for the $q$-Racah polynomials in terms of their recurrence coefficients. The explicit expression is rather complicated, especially for larger $N$, but it solves nevertheless the recurrence relation of the $q$-Racah polynomials. Without the formalism developed in this paper, the proof of this solution to the recurrence relation could be more involved. We give a simple example in Appendix \ref{sec:appendix3} to illustrate the previous remarks.  

\section{Conclusion}

In summary, we have expressed the matrix elements of some representations of the braid group on three strands in terms of $q$-Racah polynomials. To arrive at this result, we considered two intermediate Casimir elements of the tensor product of three identical spin representations of $\Usl$ and showed how they are related by a conjugation involving these braid group matrix representations. The comparison of this conjugation to the transition matrix involving $6j$-symbols allowed us to relate the braid group to the $q$-Racah polynomials. Using the invertibility of the braid matrices, we recovered the orthogonality relation of the $q$-Racah polynomials. We also obtained an addition formula for these polynomials from the braid relation, and a solution to their recurrence relation. 

Let us stress once again that the key to the results of this paper is the algebraic description of the centralizer of $\Usl$ in $\Utrp$ in terms of the universal $R$-matrix, as studied in \cite{CGVZ}. A similar work has been done for the superalgebra $\osp(1|2)$ and its universal $R$-matrix in \cite{CVZ1}. In this case, the centralizer satisfies the Bannai--Ito algebra. Therefore, we expect that the approach described in the present paper could be reproduced for this case and lead to formulas involving the Bannai--Ito polynomials. More simply, we could repeat the procedure for the algebra $U(\gsl_2)$ and the Racah polynomials, that is examine what emerges in the limit $q\to 1$. In this case, the $R$-matrix should be replaced by a permutation matrix. It could also be possible to generalize this method and obtain results involving multivariate orthogonal polynomials by considering a higher rank tensor product. We plan to examine these aspects. 

In this paper, we have mainly considered representations, since it was sufficient for our purposes. However, a more abstract point of view could have been taken by considering the algebra of the centralizer of the diagonal embedding of $\Usl$ in the tensor product of three identical spin-$s$ representations, that we denoted $\cC_s$. A conjecture that describes algebraically the centralizer of $\Usl$ in the tensor product of any three irreducible representations of spins $j_1,j_2,j_3$ was proposed in \cite{CVZ2}, where it is suggested that this more general centralizer is a quotient of the Askey--Wilson algebra. In the case $j_1=j_2=j_3=s$, the centralizer $\cC_s$ should also be a quotient of the braid group algebra on three strands; it is known to be isomorphic to the Temperley--Lieb algebra for $s=\frac{1}{2}$, and to the Birman--Murakami--Wenzl algebra for $s=1$. We trust that it should be possible to present the centralizer $\cC_s$  explicitly as a quotient of the braid group algebra, for any spin $s$, by using an algebraic version of the braid matrices $S_i$ defined in the present paper. In fact, we plan to report on this in the future. This would provide a better algebraic understanding of the connection between the braid group and the $q$-Racah polynomials.  

\vspace{5pt}

%\bigskip 

{\bf Acknowledgments:}
The authors thank Paul Terwilliger for pointing out the connections with spin Leonard pairs. N. Cramp\'e is partially supported by Agence Nationale de la
Recherche Projet AHA ANR-18-CE40-0001. The work of L. Vinet is funded in part by a discovery grant of the Natural Sciences and Engineering Research Council (NSERC) of Canada. M. Zaimi held a graduate scholarship from the Fonds de recherche du Qu\'ebec -- Nature et technologies (FRQNT) and was also partly funded by an Alexander-Graham-Bell scholarship from the NSERC.   

\appendix

\section{Some standard notations and useful identities}\label{sec:appendix1}

The $q$-shifted factorial, or $q$-Pochhammer symbol, is defined by 
\begin{equation}
	(z;q)_0:=1 \quad \text{and} \quad (z;q)_k:=\prod_{i=1}^{k}(1-zq^{i-1}) \quad \text{for } k=1,2,3,... \label{eq:qf}
\end{equation}
The following identities can be found in \cite{GR} 
\begin{align}
	&(z;q^{-1})_n=(z^{-1};q)_n(-z)^nq^{-\binom{n}{2}}, \label{eq:qfid2} \\
	&(z^{-1}q^{1-n};q)_n=(z;q)_n(-z^{-1})^nq^{-{n\choose 2}}, \label{eq:qfid3} \\
	&(z^2;q^2)_n=(z,-z;q)_n, \label{eq:qfid4} \\
	&\frac{(qz,-qz;q)_n}{(z,-z;q)_n}=\frac{(q^2z^2;q^2)_n}{(z^2;q^2)_n}=\frac{1-z^2q^{2n}}{1-z^2}. \label{eq:qfid5}
\end{align}

For $r$ and $s$ any two positive integers, the $q$-hypergeometric function ${}_{r}\phi_{s}$ is defined by 
\begin{equation}
	{}_{r}\phi_{s}\biggl(\genfrac..{0pt}{}{a_1,...,a_r}{b_1,...,b_s};q,z\biggr):=\sum_{k=0}^{\infty}\frac{(a_1,...,a_r;q)_k}{(b_1,...,b_s;q)_k}\left[ (-1)^kq^{{k\choose 2}} \right]^{1+s-r}\frac{z^r}{(q;q)_r}, \label{eq:rphis}
\end{equation}
where $(z_1,...,z_r;q)_k:=(z_1;q)_k...(z_r;q)_k$ for $k=0,1,2,...$ For $N$ a nonnegative integer, the $q$-Racah polynomials are defined by \cite{Koek}
\begin{equation}
	R_n(\mu(x);\alpha,\beta,\gamma,\delta|q):={}_{4}\phi_{3}\biggl(\genfrac..{0pt}{}{q^{-n},\al\be q^{n+1},q^{-x},\ga\de q^{x+1}}{\al q,\be\de q,\ga q};q,q\biggr), \quad n=0,1,...,N, \label{eq:qRac}
\end{equation}
where $\mu(x):=q^{-x}+\ga\de q^{x+1}$, and either $\al=q^{-N-1}$, $\be\de=q^{-N-1}$ or $\ga=q^{-N-1}$.

\section{Universal $R$-matrix of $U_q(\mathfrak{sl}_2)$ and its representations}\label{sec:appendix2}

This appendix provides more details for the proof of Proposition \ref{prop:relSiXi}.

\subsection{Universal $R$-matrix}\label{ssec:Runi}
The universal $R$-matrix of $\Usl$ is an invertible element $\cR\in \Usl \otimes \Usl$ which satisfies
\begin{equation}
	\Delta(x) \cR  = \cR \Delta^{op}(x) \quad \text{for } x\in \Usl. \label{eq:RD}
\end{equation}
In the previous equation, $\Delta^{op}$ is the opposite comultiplication and is defined by $\Delta^{op}=\tau \circ \Delta$, where $\tau(x \otimes y)=y \otimes x$ for $x,y\in \Usl$.
The universal $R$-matrix also satisfies the Yang--Baxter equation
\begin{equation}
	\cR_{12}\cR_{13}\cR_{23}=\cR_{23}\cR_{13}\cR_{12}. \label{eq:YBE}
\end{equation}
We have used the usual notations: if $\cR=\cR^\alpha\otimes \cR_\alpha$, then $\cR_{12}=\cR^\alpha\otimes \cR_\alpha \otimes 1$, $\cR_{23}=1\otimes \cR^\alpha\otimes \cR_\alpha$ and $\cR_{13}=\cR^\alpha\otimes1\otimes  \cR_\alpha$ (the sum w.r.t. $\alpha$ is understood). The universal $R$-matrix of $\Usl$ takes the following explicit form \cite{Dr}
\begin{equation}
	\cR= \sum_{n=0}^\infty \frac{(q-q^{-1})^n}{[n]_q!} q^{-n(n+1)/2} (  E \otimes F )^n  (q^{-nH}\otimes q^{nH}) q ^{2(H\otimes H)}, \label{eq:uR}
\end{equation}
where $[n]_q!:=[n]_q[n-1]_q\dots [2]_q[1]_q$ and, by convention, $[0]_q!:=1$.

\subsection{Braided $R$-matrix} \label{ssec:braidR}
We define the braided $R$-matrix in $\End(M_s^{\otimes 2})$ as
\begin{equation}
	\check R:=\pi_s^{\otimes 2}(\cR) \sigma, \label{eq:bR}
\end{equation}
where $\sigma(v\otimes w)=w \otimes v$ for any $v,w \in M_s$. 
From the Yang--Baxter equation \eqref{eq:YBE}, one can show that the braided $R$-matrix satisfies the braided Yang--Baxter equation
\begin{equation}
	\check R_1 \check R_2\check R_1= \check R_2\check R_1 \check R_2, \label{eq:bYBE}
\end{equation}
where $\check R_1=\check R \otimes \bI_{2s+1}$ and $\check R_2=\bI_{2s+1} \otimes \check R$ are elements of $\text{End}(M_s^{\otimes 3})$.

We have an explicit expression for the braided $R$-matrix. This result is known \cite{LZ}, but we give a direct proof for completeness. 

\begin{prop}
	The braided $R$-matrix is given by
	\begin{equation}
		\check R = (-1)^{2s} q^{-2s(s+1)} \sum_{r=0}^{2s} (-1)^r q^{r(r+1)} P( M_r ), \label{eq:Rc}
	\end{equation}
	where $P( M_r)$ is the projector on the representation $M_r$ of the decomposition \eqref{eq:decomp2}.
\end{prop}
\proof 
From relation \eqref{eq:RD}, one can show that $\check R$ commutes with $\pi_{s}^{\otimes 2}(\Delta(x))$ for any $x\in U_q(\gsl_2)$. 
Therefore, we must have 
\begin{equation}
\check R= \sum_{r=0}^{2s} \xi_r P( M_r ) \label{eq:Rcrep0}
\end{equation}
for some parameters $\xi_r$. Let us denote by $w_s^+$ the highest weight vector of the representation $M_s$ satisfying $Ew_s^+=0$ and $q^Hw_s^+=q^s w_s^+$. 
In the previous relations, $E$ stands for $\pi_s(E)$. We will use this abuse of notation throughout this proof.
We also introduce the vectors
\begin{equation}
\omega_r^+:=\sum_{p=0}^{2s-r}    \rho_p    F^p w_s^+ \otimes F^{2s-r-p} w_s^+ \quad \text{for } r=0,1,\dots, 2s,  \label{eq:omr}
\end{equation}
with
\begin{equation}
\rho_p:=(-1)^p q^{p(r+1)} \frac{[2s-p]_q![r+p]_q!}{[2s-r-p]_q![p]_q!}.\label{eq:rho}
\end{equation}

By a direct computation, using the identity $[E,F^n]=[n]_q F^{n-1}[2H-(n-1)]_q$, one can show that $\Delta(E)\omega_r^+=0$ and $\Delta(q^{H})\omega_r^+= q^{r}\omega_r^+$. This means that $\omega_r^+$ is the highest weight vector of the representation $M_r$ in the decomposition \eqref{eq:decomp2} of $M_s^{\otimes 2}$.
In particular, one gets $P(M_k)\omega_r^+=\delta_{kr}\omega_r^+$. Hence, the eigenvalue of $\check R$ on $\omega_r^+$ is $\xi_r$.  
Using the explicit expression \eqref{eq:uR} of the universal $R$-matrix, one gets
\begin{equation}
\check R   \omega_r^+ = 
\sum_{p=0}^{2s-r} \sum_{n=0}^{2s-r-p}\frac{\rho_p q^{2(s-r-p-n)(p-s)-n(n+1+2r)/2} (q-q^{-1})^n}{[n]_q!}  E^n F^{2s-r-p} w_s^+ \otimes F^{n+p} w_s^+ \label{eq:Rcrep1}.
\end{equation}
For any integers $n$ and $k$ such that $0\leq n \leq k \leq 2s$, one can show that
\begin{equation}
E^nF^k w_s^+ = \frac{[k]_q![2s-k+n]_q!}{[k-n]_q![2s-k]_q!}F^{k-n} w_s^+. \label{eq:EnFk}
\end{equation}
Then, by using \eqref{eq:rho} and \eqref{eq:EnFk} in \eqref{eq:Rcrep1}, one finds
\begin{align}
\check R  \omega_r^+ = 
\sum_{p=0}^{2s-r} \sum_{n=0}^{2s-r-p}A(p,n) F^{2s-r-p-n} w_s^+ \otimes F^{n+p} w_s^+ \label{eq:Rcrep2},
\end{align}
where we have defined
\begin{equation}
A(p,n):=(-1)^p\frac{(q-\qi)^n}{[n]_q!} q^{2(s-r-p-n)(p-s)-n(n+1+2r)/2+p(r+1)} \frac{[2s-p]_q![r+p+n]_q!}{[2s-r-p-n]_q![p]_q!}. \label{eq:Apn}
\end{equation}
The two sums in \eqref{eq:Rcrep2} can be rearranged as follows
\begin{equation}
\check R \omega_r^+= 
\sum_{m=0}^{2s-r} \sum_{k=0}^{2s-r-m}A(k,2s-r-m-k) F^{m} w_s^+ \otimes F^{2s-r-m} w_s^+. \label{eq:Rcrep25}
\end{equation}
Using the definitions \eqref{eq:rho} and \eqref{eq:Apn} in \eqref{eq:Rcrep25}, one gets
\begin{equation}
\check R \omega_r^+= 
q^{2s(s-r-1)+r(r+1)}\sum_{m=0}^{2s-r} (-1)^m q^{-2sm} B_m \rho_m F^{m} w_s^+ \otimes F^{2s-r-m} w_s^+, \label{eq:Rcrep3}	
\end{equation}
where
\begin{equation}
	B_m:=\sum_{k=0}^{n} (-1)^k q^{k(r+m+1)} q^{-\binom{n}{2}-\binom{k}{2}} \frac{[n]_q!}{[k]_q![n-k]_q!} \frac{[n-k +r +m]_q!}{[r+m]_q!}(q-\qi)^{n-k}
\end{equation}
and $n=2s-r-m$. In terms of $q$-shifted factorials, the coefficients $B_m$ can be rewritten as
\begin{equation}
	B_m=(-1)^nq^{-2sn}\sum_{k=0}^{n} \frac{(q^2;q^2)_n}{(q^2;q^2)_k(q^2;q^2)_{n-k}}q^{2(r+m+1)k}(q^{2(r+m+1)};q^2)_{n-k}. \label{eq:Bm}
\end{equation}
Using the $q$-analog of the binomial theorem \cite{GR}
\begin{equation}
(ab;q)_n=\sum_{k=0}^{n}\frac{(q;q)_n}{(q;q)_k(q;q)_{n-k}}b^{k}(a;q)_k(b;q)_{n-k}
\end{equation}
with the substitutions $q\to q^{2}$, $a=0$ and $b=q^{2(r+m+1)}$, equation \eqref{eq:Bm} simply reduces to
\begin{equation}
B_m=(-1)^nq^{-2sn}=(-1)^{2s-r-m}q^{-2s(2s-r-m)}.
\end{equation}
Equation \eqref{eq:Rcrep3} then becomes	
\begin{equation}
\check R \omega_r^+= (-1)^{2s} 
q^{-2s(s+1)}(-1)^{r}q^{r(r+1)}\sum_{m=0}^{2s-r} \rho_m F^{m} w_s^+ \otimes F^{2s-r-m} w_s^+ .
\end{equation}
Therefore, we can identify the eigenvalue $\xi_r=(-1)^{2s} q^{-2s(s+1)}(-1)^{r}q^{r(r+1)}$ and the result \eqref{eq:Rc} follows from \eqref{eq:Rcrep0}.
\endproof

It is directly seen from the result $\eqref{eq:Rc}$ that the matrix $S_i$ defined in \eqref{eq:Si} is proportional to the representation of $\check R_i$ in $\End(V_\ell)$, for $i=1,2$.  

\section{An explicit example for equation \eqref{eq:PX2}}\label{sec:appendix3}
In this appendix, we examine in more details equation \eqref{eq:PX2} to underscore its structure with a special case.

For $N=0$, there is only one matrix component and equation \eqref{eq:PX2} leads to a triviality. For $N=1$, we write \eqref{eq:PX2} in matrix components and use \eqref{eq:Xi} and \eqref{eq:m} to get 
\begin{equation}
	P_{ij}=-(1+q^{2a+2})q^{-4(a+1)}\ga_i\ga_j\frac{(\ga_1-1)\al_{i,j} +(\chi_{a+1}- \ga_1 \chi_{a})\de_{ij}}{\chi_{a+1}-\chi_{a}} \quad \text{for } i,j=0,1. \label{eq:Palij}
\end{equation}
If $i$ or $j$ is zero, then the $q$-Racah polynomial given in \eqref{eq:compoP} is simply one and the equation \eqref{eq:Palij} provides a relation between the recurrence coefficients $\al_{m,n}$ that can be verified. If $i=j=1$, we can use \eqref{eq:qRacparam}, \eqref{eq:betai},  \eqref{eq:aljj} and the explicit value of $\ga_r$ given in \eqref{eq:Si} to get
\begin{equation}
	R_1(\mu(1);q^{2a},q^{2a},q^{2(3a+2)},q^{-2(a+2)}|q^2)=-(1+q^{2a+2}) \frac{\al_{1,0}}{\al_{0,1}} \left( \frac{1+q^{2a+2}}{\chi_{a+1}-\chi_{a}}\al_{1,0} +1\right). \label{eq:Palbraid}
\end{equation}
As expected, equation \eqref{eq:Palbraid} expresses certain $q$-Racah polynomials in terms of their recurrence coefficients. In comparison, the repetitive application of the recurrence equation \eqref{eq:EV} leads to
\begin{equation}
	R_1(\mu(1);q^{2a},q^{2a},q^{2(3a+2)},q^{-2(a+2)}|q^2)=\frac{(\chi_{a+1}-\chi_{a})}{\al_{0,1}}+1. \label{eq:Palrec}
\end{equation}
Using the explicit expressions \eqref{eq:aljmj} and \eqref{eq:aljjm} for the recurrence coefficients and the definition \eqref{eq:casimirrep}, we find
\begin{gather}
	\al_{0,1}=-q^{-2a-3}\frac{(1-q^2)(1-q^{6(a+1)})}{(1+q^{2(a+1)})} \quad \text{and} \quad \al_{1,0}=-q^{-1}\frac{(1-q^2)(1-q^{2(a+1)})}{(1+q^{2(a+1)})} \quad \text{for } N=1, \\
	\chi_{a+1}-\chi_{a}=q^{-2a-3}(1-q^2)(1-q^{4(a+1)}).
\end{gather}
It is then straightforward to verify that \eqref{eq:Palbraid} and \eqref{eq:Palrec} reduce to the same expression, as it should.

\end{document}